\documentclass[11pt]{amsart}

\usepackage[T1]{fontenc}
\usepackage[utf8]{inputenc}
\usepackage{lmodern}
\usepackage{amsmath,amssymb,amsthm,mathtools}
\usepackage{booktabs}
\usepackage{microtype}
\usepackage[hidelinks]{hyperref}

\hypersetup{pdftitle={The Bojanov--Naidenov inequality for quartics and second derivatives}, pdfauthor={Gentian Zavalani}}

\numberwithin{equation}{section}

\newtheorem{theorem}{Theorem}
\newtheorem{lemma}[theorem]{Lemma}
\newtheorem{remark}[theorem]{Remark}

\newcommand{\R}{\mathbb{R}}

\newcommand{\Pfour}{\pi_4^\circ}
\newcommand{\norm}[1]{\left\lVert #1\right\rVert}
\newcommand{\meas}{\operatorname{meas}}

\title[The Bojanov--Naidenov inequality for quartics]{The Bojanov--Naidenov inequality for quartics and second derivatives}
\author{Gentian Zavalani}
\thanks{Technische Universit{\"a}t Dresden, Institute of Numerical Mathematics, 01062 Dresden, Germany (\url{gentian.zavalani@tu-dresden.de}).}
\date{May 29, 2026}

\begin{document}

\begin{abstract}
We settle the case $n=4$, $k=2$ of the Bojanov--Naidenov problem for algebraic polynomials.  Let $P$ be a real polynomial of degree at most four with $\norm{P}_{C[-1,1]}\leq 1$, and let $T_4(x)=8x^4-8x^2+1$.  We prove that, for every $t\geq0$,
\[
   \int_{-1}^{1} \bigl(|P''(x)|-t\bigr)_+\,dx
   \leq
   \int_{-1}^{1} \bigl(|T_4''(x)|-t\bigr)_+\,dx .
\]
This tail estimate implies
\[
   \int_{-1}^{1}\varphi(|P''(x)|)\,dx
   \leq
   \int_{-1}^{1}\varphi(|T_4''(x)|)\,dx
\]
for every nondecreasing convex function $\varphi:[0,\infty)\to\R$.  If $\varphi$ is strictly increasing and convex, equality can occur only for $P=\pm T_4$.  The proof is elementary and finite.  We interpolate at the five extremal points of $T_4$; convexity then reduces the problem to the $32$ sign choices at these nodes.  At each vertex the second derivative is a quadratic polynomial, so the remaining work is an explicit comparison of level sets.
\end{abstract}

\maketitle

\section{Introduction}

Let $\pi_n$ be the space of real algebraic polynomials of degree at most $n$, and put
\[
   \pi_n^\circ:=\{P\in\pi_n:\norm{P}_{C[-1,1]}\leq 1\}.
\]
We write
\[
   T_n(x)=\cos(n\arccos x),\qquad -1\leq x\leq 1,
\]
for the Chebyshev polynomial of the first kind.  Naidenov \cite{Naidenov2012} recorded the following problem of Bojanov.  Is it true that
\begin{equation}\label{eq:BN-general}
   \int_{-1}^{1}\varphi(|P^{(k)}(x)|)\,dx
   <
   \int_{-1}^{1}\varphi(|T_n^{(k)}(x)|)\,dx
\end{equation}
for every strictly increasing convex function $\varphi$ on $[0,\infty)$, every $P\in\pi_n^\circ\setminus\{\pm T_n\}$, and every $k=2,\ldots,n-1$?

The endpoint cases have a different nature.  The case $k=n$ follows from a classical extremal property of the Chebyshev polynomial, while the case $k=1$ was proved by Bojanov \cite{Bojanov1983}; related Markov-type results had appeared earlier in \cite{Bojanov1982longest,Bojanov1982}.  The problem is also known in some special classes.  Bojanov and Rahman \cite{BojanovRahman1995} treated polynomials whose zeros lie in $[-1,1]$, and Avvakumova \cite{Avvakumova2000} proved the case $k=2$ for a larger intermediate class.  More recent work studies Bojanov--Naidenov-type questions in other settings, including differentiable functions on the real line, trigonometric polynomials, and splines, together with links to sharp Kolmogorov- and Bernstein-type inequalities \cite{Kofanov2023,Kofanov2024,Kofanov2025}.  These papers deal with different spaces or metrics from \eqref{eq:BN-general}.  The question considered here is the original algebraic-polynomial problem on $[-1,1]$.

This paper proves the full quartic case, namely $n=4$ and $k=2$.  In fact, we prove a little more than the convex integral inequality: the integrated tails of $|P''|$ are dominated by those of $|T_4''|$.  This is an increasing-convex-order statement, and it gives \eqref{eq:BN-general} in the quartic case.

The argument is quite specific to degree four.  A quartic is determined by its values at the five extremal points of $T_4$.  After this interpolation, a simple convexity argument allows us to enlarge the feasible set to the full cube $[-1,1]^5$ and then pass to its vertices.  There are only $32$ of them, and their second derivatives reduce, up to sign and reflection, to a short list of quadratic polynomials.  The rest of the proof compares the level sets of these quadratics with the level sets of $T_4''$.  Thus the finite character of the proof is also its limitation: it gives a complete solution of this concrete intermediate case, but it is not a direct argument for higher degrees.

\begin{theorem}\label{thm:main}
Let $P\in\pi_4$ satisfy $\norm{P}_{C[-1,1]}\leq 1$, and let
\[
   T_4(x)=8x^4-8x^2+1.
\]
Then, for every $t\geq 0$,
\begin{equation}\label{eq:tail-main}
   \int_{-1}^{1}\bigl(|P''(x)|-t\bigr)_+\,dx
   \leq
   \int_{-1}^{1}\bigl(|T_4''(x)|-t\bigr)_+\,dx .
\end{equation}
Consequently, for every nondecreasing convex function $\varphi:[0,\infty)\to\R$,
\begin{equation}\label{eq:convex-main}
   \int_{-1}^{1}\varphi(|P''(x)|)\,dx
   \leq
   \int_{-1}^{1}\varphi(|T_4''(x)|)\,dx .
\end{equation}
If $P\neq\pm T_4$ and $\varphi$ is strictly increasing and convex, then \eqref{eq:convex-main} is strict.
\end{theorem}

We briefly outline the proof.  For a fixed tail parameter $t$, the map from the five interpolating values to the integrated tail of $|P''|$ is convex.  Its maximum on $[-1,1]^5$ is therefore attained at a vertex.  At the vertices we only have to handle finitely many quadratic polynomials, up to sign and reflection; these are compared with
\[
   T_4''(x)=96x^2-16.
\]

\section{A tail-to-convex lemma}

We start with a standard tail representation for increasing convex functionals.  The lemma is stated with the strict part included, because this will be needed at the end of the proof.

\begin{lemma}\label{lem:tail-to-convex}
Let $f,g$ be measurable functions on a finite measure space, with values in $[0,M]$.  Suppose
\begin{equation}\label{eq:tail-order-general}
   \int (f-t)_+\,d\nu\leq \int (g-t)_+\,d\nu,
   \qquad 0\leq t\leq M.
\end{equation}
Then
\[
   \int \varphi(f)\,d\nu\leq \int \varphi(g)\,d\nu
\]
for every nondecreasing convex function $\varphi:[0,M]\to\R$.
If, in addition,
\begin{equation}\label{eq:strict-at-zero}
   \int f\,d\nu<\int g\,d\nu,
\end{equation}
then the last inequality is strict for every strictly increasing convex function $\varphi$.
\end{lemma}

\begin{proof}
Every convex function on $[0,M]$ can be written in the form
\begin{equation}\label{eq:hinge-representation}
   \varphi(u)=\varphi(0)+mu+
   \int_{[0,M]}(u-s)_+\,d\mu(s),
\end{equation}
where $m\in\R$ and $\mu$ is a nonnegative finite Borel measure.  When $\varphi$ is nondecreasing, we may choose the linear coefficient with $m\geq0$, after absorbing a possible atom at $0$ into $m$.  Integrating \eqref{eq:hinge-representation} and applying \eqref{eq:tail-order-general} proves the non-strict inequality.

Now assume \eqref{eq:strict-at-zero}.  Put
\[
   D(t):=\int(g-t)_+\,d\nu-\int(f-t)_+\,d\nu .
\]
Then $D$ is continuous and $D(0)>0$.  Hence $D(t)>0$ on some interval $0\leq t\leq\delta$.  If $m>0$, the strict inequality already follows from the linear term.  If $m=0$ and $\varphi$ is strictly increasing, then $\mu((0,\delta))>0$; otherwise the representation \eqref{eq:hinge-representation} would make $\varphi$ constant on $[0,\delta]$.  Therefore
\[
   \int \varphi(g)\,d\nu-\int \varphi(f)\,d\nu
   = mD(0)+\int_{[0,M]}D(s)\,d\mu(s)>0.
\]
\end{proof}

\section{Reduction to vertices}

Let
\[
   x_j=\cos\frac{j\pi}{4},\qquad j=0,1,2,3,4.
\]
Thus
\[
   x_0=1,\qquad x_1=\frac{\sqrt2}{2},\qquad x_2=0,
   \qquad x_3=-\frac{\sqrt2}{2},\qquad x_4=-1.
\]
For $y=(y_0,\ldots,y_4)\in\R^5$, let $P_y$ be the unique polynomial of degree at most four satisfying
\[
   P_y(x_j)=y_j,\qquad j=0,\ldots,4.
\]
For fixed $t\geq 0$, define
\begin{equation}\label{eq:Ft}
   F_t(y):=\int_{-1}^{1}\bigl(|P_y''(x)|-t\bigr)_+\,dx .
\end{equation}

\begin{lemma}\label{lem:convex-reduction}
For each $t\geq 0$, the map $F_t:\R^5\to\R$ is convex.  Consequently,
\[
   \max_{y\in[-1,1]^5}F_t(y)
\]
is attained at a vertex of the cube $[-1,1]^5$.
\end{lemma}

\begin{proof}
Interpolation depends linearly on the nodal values.  Thus, for fixed $x$, the map $y\mapsto P_y''(x)$ is linear.  The scalar function
\[
   u\mapsto (|u|-t)_+=\max\{u-t,-u-t,0\}
\]
is convex, so $y\mapsto (|P_y''(x)|-t)_+$ is convex as well.  Integrating over $[-1,1]$ keeps convexity, which proves the first assertion.

For the second assertion, write an arbitrary point of the cube as a convex combination of vertices, say $y=\sum_v\lambda_v v$.  Then
\[
   F_t(y)\leq \sum_v\lambda_vF_t(v)\leq \max_vF_t(v).
\]
Hence a maximum is attained at some vertex.
\end{proof}

The nodal vectors coming from polynomials in $\Pfour$ form only a subset of $[-1,1]^5$.  We shall nevertheless work on the full cube.  This enlargement is harmless, since we only need an upper bound.

\section{The vertex second derivatives}

Let $\ell_j$ be the Lagrange basis polynomial associated with $x_j$.  A direct calculation gives
\begin{align*}
   \ell_0''(x)&=12x^2+6x-1,\\
   \ell_1''(x)&=-24x^2-6\sqrt2\,x+4,\\
   \ell_2''(x)&=24x^2-6,\\
   \ell_3''(x)&=-24x^2+6\sqrt2\,x+4,\\
   \ell_4''(x)&=12x^2-6x-1.
\end{align*}
Therefore
\[
   P_y''(x)=\sum_{j=0}^{4}y_j\ell_j''(x)=A(y)x^2+B(y)x+C(y),
\]
where
\begin{align}
   A(y)&=12y_0-24y_1+24y_2-24y_3+12y_4,\label{eq:A}\\
   B(y)&=6y_0-6\sqrt2\,y_1+6\sqrt2\,y_3-6y_4,\label{eq:B}\\
   C(y)&=-y_0+4y_1-6y_2+4y_3-y_4.\label{eq:C}
\end{align}

At the vertices we have $y_j\in\{-1,1\}$, hence $32$ possible second derivatives.  The transformations $q(x)\mapsto -q(x)$ and $q(x)\mapsto q(-x)$ leave all level-set and tail quantities unchanged.  After using these symmetries and removing duplicates, the list is:
\begin{equation}\label{eq:vertex-list}
\begin{array}{c|l}
\text{type} & q(x)=P_y''(x) \\ \hline
0 & 0 \\
L & 12\sqrt2\,x+4 \\
Q_1 & 24x^2+12x-2 \\
Q_2 & 48x^2+12\sqrt2\,x-8 \\
Q_3 & 24x^2+12(\sqrt2-1)x-6 \\
Q_4 & 48x^2-12 \\
Q_5 & 72x^2+12x-14 \\
Q_6 & 24x^2+12(\sqrt2+1)x-6 \\
Q_7 & 48x^2-4 \\
Q_* & 96x^2-16.
\end{array}
\end{equation}
The last row is $T_4''$ itself.  It comes from the Chebyshev nodal vector
\[
   (1,-1,1,-1,1),
\]
and its negative corresponds to $-T_4$.

We shall need the following elementary data.  For a quadratic row we write
\[
   q(x)=a(x-h)^2-v.
\]

\begin{table}[htbp]
\caption{Data for the vertex derivatives.  The norm is the uniform norm on $[-1,1]$.}
\label{tab:vertex-data}
\centering
\footnotesize
\begin{tabular}{@{}ccccccc@{}}
\toprule
Type & $a$ & $h$ & $v$ & $q(-1)$ & $q(1)$ & $\norm{q}_{C[-1,1]}$ \\
\midrule
$0$ & -- & -- & -- & $0$ & $0$ & $0$ \\
$L$ & -- & -- & -- & $4-12\sqrt2$ & $4+12\sqrt2$ & $4+12\sqrt2$ \\
$Q_1$ & $24$ & $-\frac14$ & $\frac72$ & $10$ & $34$ & $34$ \\
$Q_2$ & $48$ & $-\frac{\sqrt2}{8}$ & $\frac{19}{2}$ & $40-12\sqrt2$ & $40+12\sqrt2$ & $40+12\sqrt2$ \\
$Q_3$ & $24$ & $\frac{1-\sqrt2}{4}$ & $\frac{21}{2}-3\sqrt2$ & $30-12\sqrt2$ & $6+12\sqrt2$ & $6+12\sqrt2$ \\
$Q_4$ & $48$ & $0$ & $12$ & $36$ & $36$ & $36$ \\
$Q_5$ & $72$ & $-\frac1{12}$ & $\frac{29}{2}$ & $46$ & $70$ & $70$ \\
$Q_6$ & $24$ & $-\frac{1+\sqrt2}{4}$ & $\frac{21}{2}+3\sqrt2$ & $6-12\sqrt2$ & $30+12\sqrt2$ & $30+12\sqrt2$ \\
$Q_7$ & $48$ & $0$ & $4$ & $44$ & $44$ & $44$ \\
$Q_*$ & $96$ & $0$ & $16$ & $80$ & $80$ & $80$ \\
\bottomrule
\end{tabular}
\end{table}

\begin{lemma}\label{lem:markov-quartic}
If $P\in\Pfour$, then
\[
   \norm{P''}_{C[-1,1]}\leq 80=\norm{T_4''}_{C[-1,1]}.
\]
Equality is possible only for $P=\pm T_4$.
\end{lemma}

\begin{proof}
Let $y_j=P(x_j)$.  Then $y\in[-1,1]^5$ and $P=P_y$.  The map
\[
   y\mapsto \norm{P_y''}_{C[-1,1]}
\]
is convex, since it is the supremum, over $x\in[-1,1]$, of the convex functions $y\mapsto |P_y''(x)|$.  Its maximum on the cube is therefore attained at a vertex.  Table \ref{tab:vertex-data} shows that the two Chebyshev vertices have norm $80$, while every other vertex has norm strictly smaller than $80$.

Suppose equality holds for some $y\in[-1,1]^5$.  Write $y=\sum_v\lambda_v v$ as a convex combination of vertices.  Convexity gives
\[
   80=\norm{P_y''}_{C[-1,1]}
   \leq \sum_v\lambda_v\norm{P_v''}_{C[-1,1]}
   \leq 80.
\]
Thus every vertex with $\lambda_v>0$ must be a Chebyshev vertex; any other vertex would make the last inequality strict.  Hence $y$ lies on the segment joining $(1,-1,1,-1,1)$ to $(-1,1,-1,1,-1)$.  In other words, $P=\alpha T_4$ with $|\alpha|\leq1$, and equality in the norm forces $|\alpha|=1$.
\end{proof}

\section{The Chebyshev tail}

Set
\[
   q_*(x):=T_4''(x)=96x^2-16.
\]
For a measurable function $q$ on $[-1,1]$, define the level-set function
\[
   \Lambda_q(t):=\meas\{x\in[-1,1]: |q(x)|>t\}.
\]
Then
\begin{equation}\label{eq:tail-level}
   \int_{-1}^{1}(|q(x)|-t)_+\,dx
   =\int_t^\infty \Lambda_q(s)\,ds .
\end{equation}
We write $\Lambda_*:=\Lambda_{q_*}$.  The polynomial $q_*$ vanishes at $|x|=1/\sqrt6$, takes the value $-16$ at $x=0$, and takes the value $80$ at $x=\pm1$.  Hence
\begin{equation}\label{eq:Lambda-star}
\Lambda_*(t)=
\begin{cases}
\displaystyle
2-2\sqrt{\frac{16+t}{96}}
+2\sqrt{\frac{16-t}{96}},
&0\leq t\leq16,\\[1.1em]
\displaystyle
2-2\sqrt{\frac{16+t}{96}},
&16\leq t\leq80,\\[1.1em]
0,
&t\geq80.
\end{cases}
\end{equation}
Equivalently, the integrated Chebyshev tail is
\begin{equation}\label{eq:R-star}
\int_{-1}^{1}(|q_*(x)|-t)_+\,dx=
\begin{cases}
\displaystyle
32-2t+
\frac{(16+t)^{3/2}+(16-t)^{3/2}}{3\sqrt6},
&0\leq t\leq16,\\[1.1em]
\displaystyle
32-2t+
\frac{(16+t)^{3/2}}{3\sqrt6},
&16\leq t\leq80,\\[1.1em]
0,
&t\geq80.
\end{cases}
\end{equation}

\section{Finite level-set comparison}

\begin{lemma}\label{lem:standard-quadratic-levels}
Let
\[
   q(x)=a(x-h)^2-v,
   \qquad a>0,
   \qquad -1<h\leq0,
   \qquad v>0,
\]
and set
\[
   L:=q(-1),\qquad R:=q(1).
\]
Assume
\[
   0<v\leq L\leq R.
\]
Then
\[
\Lambda_q(t)=
\begin{cases}
\displaystyle
2-2\sqrt{\frac{v+t}{a}}+2\sqrt{\frac{v-t}{a}},
&0\leq t\leq v,\\[1.1em]
\displaystyle
2-2\sqrt{\frac{v+t}{a}},
&v\leq t\leq L,\\[1.1em]
\displaystyle
1-h-\sqrt{\frac{v+t}{a}},
&L\leq t\leq R,\\[1.1em]
0,&t\geq R.
\end{cases}
\]
\end{lemma}

\begin{proof}
The minimum of $q$ is $-v$, attained at $x=h$.  If $0\leq t\leq v$, then the set where $q<-t$ is
\[
   h-\sqrt{\frac{v-t}{a}}<x<h+\sqrt{\frac{v-t}{a}},
\]
and the set where $q>t$ is the union of the two outer intervals cut off by
\[
   x=h\pm\sqrt{\frac{v+t}{a}}.
\]
The assumptions $0<v\leq L\leq R$ ensure that the relevant roots lie in $[-1,1]$ for this range of $t$.  Adding the lengths gives the first formula.

If $v\leq t\leq L$, the negative interval has disappeared, and both roots of $q=t$ still lie in $[-1,1]$.  This gives the second formula.  If $L\leq t\leq R$, the left root of $q=t$ has moved past $-1$, while the right root remains in $[-1,1]$ until $t=R$.  Thus only the right outer interval remains, which gives the third formula.  For $t\geq R$ the set is empty.
\end{proof}

The remaining point is the following finite comparison.

\begin{lemma}\label{lem:finite-level}
For each non-Chebyshev vertex derivative $q$ in \eqref{eq:vertex-list},
\begin{equation}\label{eq:level-dom}
   \Lambda_q(t)\leq \Lambda_*(t),\qquad t\geq0.
\end{equation}
Consequently,
\begin{equation}\label{eq:tail-dom-vertices}
   \int_{-1}^{1}(|q(x)|-t)_+\,dx
   \leq
   \int_{-1}^{1}(|T_4''(x)|-t)_+\,dx,
   \qquad t\geq0.
\end{equation}
\end{lemma}

\begin{proof}
The zero vertex is immediate, so we consider the remaining cases.

\subsection*{The linear vertex}
Let
\[
   q(x)=12\sqrt2\,x+4,
   \qquad m:=12\sqrt2.
\]
Then
\[
\Lambda_q(t)=
\begin{cases}
\displaystyle 2-\frac{2t}{m},
&0\leq t\leq m-4,\\[0.9em]
\displaystyle \frac{m+4-t}{m},
&m-4\leq t\leq m+4,\\[0.9em]
0,&t\geq m+4.
\end{cases}
\]
For $0\leq t\leq m-4$, we need to prove
\[
   \frac{t}{m}
   \geq
   \sqrt{\frac{16+t}{96}}-
   \sqrt{\frac{16-t}{96}}.
\]
For $t>0$ the difference of square roots on the right is
\[
   \frac{2t}{\sqrt{96}\bigl(\sqrt{16+t}+\sqrt{16-t}\bigr)}.
\]
It is therefore enough that
\[
   \sqrt{96}\bigl(\sqrt{16+t}+\sqrt{16-t}\bigr)
   \geq 2m,
\]
which holds because the left-hand side is at least $16\sqrt6$, while $2m=24\sqrt2$.  The endpoint $t=0$ follows by continuity.

For $m-4\leq t\leq m+4$, we have $\Lambda_q(t)\leq 8/m=\sqrt2/3$.  Since $\Lambda_*$ is decreasing, it remains only to check
\[
   \frac{\sqrt2}{3}\leq \Lambda_*(m+4)
   =2-2\sqrt{\frac{20+12\sqrt2}{96}}.
\]
Both sides are nonnegative; after squaring, this is equivalent to $61>33\sqrt2$.  Hence \eqref{eq:level-dom} holds for the linear vertex.

\subsection*{Six standard quadratic vertices}
Consider the six rows
\[
   Q_1,\quad Q_2,
   \quad Q_3,
   \quad Q_4,
   \quad Q_5,
   \quad Q_7.
\]
Write
\[
   q(x)=a(x-h)^2-v,
   \qquad a>0,
   \qquad v>0,
   \qquad h\leq0,
\]
and set
\[
   L:=q(-1),\qquad R:=q(1),\qquad L\leq R.
\]
For these six rows, Table \ref{tab:vertex-data} gives the inequalities
\begin{equation}\label{eq:standard-conditions}
\begin{gathered}
   0<v\leq16,\qquad a\leq96,\qquad -1<h\leq0,\\
   a\leq3(v+16),\qquad 0<v\leq L\leq R<80.
\end{gathered}
\end{equation}
Thus Lemma \ref{lem:standard-quadratic-levels} applies to each of the six rows.  We compare the level sets on the intervals $[0,v]$, $[v,L]$, $[L,R]$, and $[R,\infty)$.

If $0\leq t\leq v$, then Lemma \ref{lem:standard-quadratic-levels} gives
\[
   \Lambda_q(t)=2-2\sqrt{\frac{v+t}{a}}+2\sqrt{\frac{v-t}{a}}.
\]
We have to show
\[
\sqrt{\frac{v+t}{a}}-
\sqrt{\frac{v-t}{a}}
\geq
\sqrt{\frac{16+t}{96}}-
\sqrt{\frac{16-t}{96}}.
\]
The function $s\mapsto\sqrt{s+t}-\sqrt{s-t}$ is decreasing for $s\geq t$.  Since $v\leq16$ and $a\leq96$, the required inequality follows.

If $v\leq t\leq L$, then Lemma \ref{lem:standard-quadratic-levels} gives
\[
   \Lambda_q(t)=2-2\sqrt{\frac{v+t}{a}}.
\]
First assume $t\leq16$.  It remains to prove
\[
\sqrt{\frac{v+t}{a}}
\geq
\sqrt{\frac{16+t}{96}}-
\sqrt{\frac{16-t}{96}}.
\]
Using
\[
\left(
\sqrt{\frac{16+t}{96}}-
\sqrt{\frac{16-t}{96}}
\right)^2
\leq \frac{t^2}{768},
\]
it is enough to show
\[
   \frac{v+t}{a}\geq \frac{t^2}{768}.
\]
The function $t^2/(v+t)$ is increasing for $t\geq0$, and $t\leq16$, so
\[
   \frac{t^2}{v+t}\leq \frac{256}{v+16}.
\]
Thus the claim follows from $a\leq3(v+16)$.

If $16\leq t\leq L$, then
\[
   \Lambda_*(t)=2-2\sqrt{\frac{16+t}{96}},
\]
and we only have to prove
\[
   \frac{v+t}{a}\geq \frac{16+t}{96}.
\]
Equivalently,
\[
   96(v+t)\geq a(16+t).
\]
Because $a\leq96$, the difference between the left and right sides is increasing in $t$.  It is therefore enough to check $t=16$, where the condition becomes $a\leq3(v+16)$.

It remains to consider $L\leq t\leq R$.  In this interval, Lemma \ref{lem:standard-quadratic-levels} gives
\[
   \Lambda_q(t)=1-h-\sqrt{\frac{v+t}{a}}.
\]
For the two rows with $L<16$, namely
\[
   q(x)=24x^2+12x-2,
   \qquad
   q(x)=24x^2+12(\sqrt2-1)x-6,
\]
we have $\Lambda_q(t)\leq -2h$.  Since $\Lambda_*$ is decreasing, it remains to verify $-2h\leq\Lambda_*(R)$.

For $q(x)=24x^2+12x-2$, this is
\[
   \frac12\leq 2-2\sqrt{\frac{50}{96}}.
\]
Both sides are nonnegative, and squaring reduces this to $25\leq27$.  For $q(x)=24x^2+12(\sqrt2-1)x-6$, the desired inequality is
\[
   \frac{\sqrt2-1}{2}
   \leq
   2-2\sqrt{\frac{22+12\sqrt2}{96}}.
\]
Both sides are nonnegative, and after squaring this becomes $72\sqrt2\leq140$.

For the two rows with $16<L<R$, namely
\[
   q(x)=48x^2+12\sqrt2\,x-8,
   \qquad
   q(x)=72x^2+12x-14,
\]
put $E(t):=\Lambda_*(t)-\Lambda_q(t)$.  On $[L,R]$, one has $t\geq16$, and $E'(t)<0$ provided
\[
   a(v+t)>24(16+t).
\]
This condition need only be checked at $t=L$.  For the two rows, the left side minus the right side equals, respectively,
\[
   1032-288\sqrt2>0,
   \qquad
   2868>0.
\]
Hence $E$ is decreasing on $[L,R]$, and so
\[
   E(t)\geq E(R)=\Lambda_*(R)>0,
\]
because $R<80$ and $\Lambda_q(R)=0$.  The remaining two standard rows, $48x^2-12$ and $48x^2-4$, have $L=R$, so there is no endpoint interval to check.  This proves the comparison for the six standard quadratic vertices.

\subsection*{The exceptional quadratic vertex}
It remains to handle the row
\[
   q(x)=24x^2+12(\sqrt2+1)x-6.
\]
Write
\[
   q(x)=24(x-h)^2-v,
   \qquad
   h=-\frac{1+\sqrt2}{4},
   \qquad
   v=\frac{21}{2}+3\sqrt2.
\]
Then
\[
   q(-1)=6-12\sqrt2<0,
   \qquad
   q(1)=30+12\sqrt2<80.
\]
Set
\[
   s_0:=12\sqrt2-6.
\]
We have
\[
   0<s_0<v<16<30+12\sqrt2<80.
\]

For $0\leq t\leq s_0$,
\[
   \Lambda_q(t)=2+
   \sqrt{\frac{v-t}{24}}-
   \sqrt{\frac{v+t}{24}}.
\]
Since $v<16$ and $s\mapsto\sqrt{s+t}-\sqrt{s-t}$ is decreasing in $s$, we have
\[
   \sqrt{v+t}-\sqrt{v-t}
   \geq
   \sqrt{16+t}-\sqrt{16-t}.
\]
Dividing by $\sqrt{24}$ gives
\[
\sqrt{\frac{v+t}{24}}-
\sqrt{\frac{v-t}{24}}
\geq
2\left(
\sqrt{\frac{16+t}{96}}-
\sqrt{\frac{16-t}{96}}
\right),
\]
and hence $\Lambda_q(t)\leq\Lambda_*(t)$ on $[0,s_0]$.

For $s_0\leq t\leq v$, the negative part is fully inside $[-1,1]$, and
\[
   \Lambda_q(t)=1-h-
   \sqrt{\frac{v+t}{24}}+
   2\sqrt{\frac{v-t}{24}}.
\]
Let $E(t):=\Lambda_*(t)-\Lambda_q(t)$.  Then
\[
E'(t)=\frac1{\sqrt{96}}
\left(
-\frac1{\sqrt{16+t}}
-\frac1{\sqrt{16-t}}
+\frac1{\sqrt{v+t}}
+\frac2{\sqrt{v-t}}
\right)>0,
\]
because $v<16$.  Since the comparison has already been proved at $t=s_0$, it holds on $[s_0,v]$.

For $v\leq t\leq16$,
\[
   \Lambda_q(t)=1-h-\sqrt{\frac{v+t}{24}}.
\]
Again set $E(t):=\Lambda_*(t)-\Lambda_q(t)$.  Then
\[
E'(t)=\frac1{\sqrt{96}}
\left(
-\frac1{\sqrt{16+t}}
-\frac1{\sqrt{16-t}}
+\frac1{\sqrt{v+t}}
\right)<0,
\]
since $t\geq v$ implies $v+t>16-t$.  Thus $E$ is decreasing on $[v,16]$, so it remains to check $E(16)>0$.  At $t=16$,
\[
E(16)=
2-\frac2{\sqrt3}
-\frac{5+\sqrt2}{4}
+
\sqrt{\frac{53+6\sqrt2}{48}}.
\]
Using
\[
   \sqrt{\frac{53+6\sqrt2}{48}}>\frac{11}{10},
   \qquad
   \sqrt2<\frac32,
   \qquad
   \frac2{\sqrt3}<\frac76,
\]
we get
\[
   E(16)>
   2-\frac76-\frac54-\frac38+\frac{11}{10}
   =\frac{37}{120}>0.
\]
Hence the comparison holds on $[v,16]$.

For $16\leq t\leq30+12\sqrt2$,
\[
   \Lambda_q(t)=1-h-\sqrt{\frac{v+t}{24}},
   \qquad
   \Lambda_*(t)=2-2\sqrt{\frac{16+t}{96}}.
\]
Thus
\[
\frac{d}{dt}(\Lambda_*(t)-\Lambda_q(t))
=
\frac1{\sqrt{96}}
\left(
\frac1{\sqrt{v+t}}-
\frac1{\sqrt{16+t}}
\right)>0,
\]
because $v<16$.  Since the difference is positive at $t=16$, it remains positive up to $30+12\sqrt2$.  For larger $t$, $\Lambda_q(t)=0$, so the comparison is trivial.  This completes the proof of \eqref{eq:level-dom}.

Finally, \eqref{eq:tail-dom-vertices} follows from \eqref{eq:tail-level}.
\end{proof}

\begin{lemma}\label{lem:L1-strict}
If $P\in\Pfour$ and $P\neq\pm T_4$, then
\begin{equation}\label{eq:L1-strict}
   \int_{-1}^{1}|P''(x)|\,dx
   <
   \int_{-1}^{1}|T_4''(x)|\,dx .
\end{equation}
\end{lemma}

\begin{proof}
Let $q$ be a non-Chebyshev vertex derivative.  Lemma \ref{lem:finite-level} gives $\Lambda_q\leq\Lambda_*$, and Table \ref{tab:vertex-data} gives $\norm{q}_{C[-1,1]}<80$.  Since $\Lambda_*(s)>0$ for $0\leq s<80$, formula \eqref{eq:tail-level} with $t=0$ gives
\[
   \int_{-1}^{1}|q(x)|\,dx
   <
   \int_{-1}^{1}|T_4''(x)|\,dx
\]
for every non-Chebyshev vertex.  The two Chebyshev vertices give equality.

The functional $F_0(y)=\int_{-1}^{1}|P_y''(x)|\,dx$ is convex.  If equality held in \eqref{eq:L1-strict}, then in a convex decomposition of the nodal vector $y=(P(x_0),\ldots,P(x_4))$ no non-Chebyshev vertex could appear with positive weight; otherwise the convexity estimate would be strict.  Thus $y$ would lie on the segment between $(1,-1,1,-1,1)$ and its negative.  Hence $P=\alpha T_4$ with $|\alpha|\leq1$, and equality in \eqref{eq:L1-strict} would force $|\alpha|=1$.  This contradicts $P\neq\pm T_4$, and the lemma follows.
\end{proof}

\section{Proof of the main theorem}

Fix $t\geq0$.  Lemma \ref{lem:finite-level}, together with the row $q=T_4''$, shows that every vertex $v\in\{-1,1\}^5$ satisfies
\[
   F_t(v)
   \leq
   \int_{-1}^{1}\bigl(|T_4''(x)|-t\bigr)_+\,dx .
\]
By Lemma \ref{lem:convex-reduction}, the same estimate holds for all $y\in[-1,1]^5$.  Now let $P\in\Pfour$ and take $y_j=P(x_j)$.  Then $y\in[-1,1]^5$ and $P=P_y$, so
\[
   \int_{-1}^{1}\bigl(|P''(x)|-t\bigr)_+\,dx
   \leq
   \int_{-1}^{1}\bigl(|T_4''(x)|-t\bigr)_+\,dx.
\]
This proves \eqref{eq:tail-main}.

By Lemma \ref{lem:markov-quartic}, both $|P''|$ and $|T_4''|$ take their values in $[0,80]$.  We may therefore apply Lemma \ref{lem:tail-to-convex} with
\[
   f(x)=|P''(x)|,
   \qquad
   g(x)=|T_4''(x)|.
\]
This gives \eqref{eq:convex-main} for every nondecreasing convex $\varphi$.

Finally suppose that $P\neq\pm T_4$.  Lemma \ref{lem:L1-strict} gives
\[
   \int_{-1}^{1}|P''(x)|\,dx
   <
   \int_{-1}^{1}|T_4''(x)|\,dx .
\]
The strict part of Lemma \ref{lem:tail-to-convex} then gives strict inequality in \eqref{eq:convex-main} for every strictly increasing convex $\varphi$.  This completes the proof.

\begin{remark}
The equality cases for the whole tail family are also exactly $P=\pm T_4$.  Indeed, if equality held in \eqref{eq:tail-main} for all $t$, then the highest tail levels would force $\norm{P''}_{C[-1,1]}=80$.  Lemma \ref{lem:markov-quartic} would then give $P=\pm T_4$.
\end{remark}

\appendix

\section{Audit of the vertex enumeration}
\label{app:vertex-enumeration}
This appendix gives the finite enumeration used in the reduction to \eqref{eq:vertex-list}.  Reflection of the interval sends a nodal vector
\[
   (y_0,y_1,y_2,y_3,y_4)
   \quad\text{to}\quad
   (y_4,y_3,y_2,y_1,y_0),
\]
and multiplication of the polynomial by $-1$ sends $y$ to $-y$.  Both operations preserve the level-set and integrated-tail quantities, since they replace $q(x)$ by $q(-x)$ or by $-q(x)$.  The representatives below, together with these symmetries, account for all $32$ vertices; the orbit sizes in the third column sum to $32$.

\begin{table}[htbp]
\caption{Representative nodal vectors for the vertex derivatives.}
\label{tab:vertex-enumeration}
\centering
\footnotesize
\begin{tabular}{@{}ccll@{}}
\toprule
Type & representative $y$ & orbit size & $P_y''(x)$ \\
\midrule
$0$ & $(-1,-1,-1,-1,-1)$ & $2$ & $0$ \\
$L$ & $(1,-1,-1,1,1)$ & $4$ & $12\sqrt2\,x+4$ \\
$Q_1$ & $(1,-1,-1,-1,-1)$ & $4$ & $24x^2+12x-2$ \\
$Q_2$ & $(1,-1,1,1,1)$ & $4$ & $48x^2+12\sqrt2\,x-8$ \\
$Q_3$ & $(-1,-1,1,1,1)$ & $4$ & $24x^2+12(\sqrt2-1)x-6$ \\
$Q_4$ & $(-1,-1,1,-1,-1)$ & $2$ & $48x^2-12$ \\
$Q_5$ & $(1,-1,1,-1,-1)$ & $4$ & $72x^2+12x-14$ \\
$Q_6$ & $(1,-1,1,1,-1)$ & $4$ & $24x^2+12(\sqrt2+1)x-6$ \\
$Q_7$ & $(1,-1,-1,-1,1)$ & $2$ & $48x^2-4$ \\
$Q_*$ & $(1,-1,1,-1,1)$ & $2$ & $96x^2-16$ \\
\bottomrule
\end{tabular}
\end{table}

Each row is obtained by substituting the representative vector into \eqref{eq:A}--\eqref{eq:C}.  Thus Table \ref{tab:vertex-enumeration} gives a direct check of the reduction from the $32$ vertex derivatives to the list \eqref{eq:vertex-list}.

\end{document}